\documentclass[10pt]{article}
\usepackage{amsfonts}
\usepackage{subfigure} 
\usepackage{epsfig}
\usepackage{amsmath, amsthm}
\newcommand{\be}{\begin{eqnarray}}     	\newcommand{\ee}{\end{eqnarray}}

\newcommand{\vol}{\mathrm{Vol}}
\newcommand{\dist}{\mathrm{dist}}
\newcommand{\rem}{\mathrm{Rm}}
\newcommand{\ric}{\mathrm{Ric}}
\newcommand{\diam}{\mathrm{diam}}

\newcommand{\supp}{\mathrm{supp}}
\newcommand{\const}{\mathrm{const}}
 
\title{The compactness result for K\"ahler Ricci solitons}
\author{Huai Dong Cao, N.Sesum}

\date{} 
\theoremstyle{plain} 
\newtheorem{dummy}{Dummy}

\theoremstyle{definition}
\newtheorem{definition}[dummy]{Definition}

\theoremstyle{plain}
\newtheorem{corollary}[dummy]{Corollary}
\newtheorem{remark}[dummy]{Remark}
\newtheorem{lemma}[dummy]{Lemma}
\newtheorem{theorem}[dummy]{Theorem}
\newtheorem{proposition}[dummy]{Proposition}

\begin{document}

\maketitle

\begin{abstract}
In this paper we prove the compactness result for compact K\"ahler
Ricci gradient shrinking solitons. If $(M_i,g_i)$ is a sequence of
K\"ahler Ricci solitons of real dimension $n \ge 4$, whose curvatures
have uniformly bounded $L^{n/2}$ norms, whose Ricci curvatures are
uniformly bounded from below and $\mu(g_i,1/2) \ge A$ (where $\mu$ is
Perelman's functional), there is a subsequence $(M_i,g_i)$ converging
to a compact orbifold $(M_{\infty},g_{\infty})$ with finitely many
isolated singularitites, where $g_{\infty}$ is a K\"ahler Ricci
soliton metric in an orbifold sense (satisfies a soliton equation away
from singular points and smoothly extends in some gauge to a metric
satisfying K\"ahler Ricci soliton equation in a lifting around
singular points).  
\end{abstract}

\begin{section}{Introduction}

Let $g(t)$ be a K\"ahler Ricci flow on a compact, K\"ahler manifold
$M$,
\begin{equation}
\label{equation-kr_flow}
\frac{d}{dt}g_{i\bar{j}} = g_{i\bar{j}} - R_{i\bar{j}} =
u_{i\bar{j}}.
\end{equation}
Very special solutions of (\ref{equation-kr_flow}) are K\"ahler Ricci
solitons. 

\begin{definition}
A K\"ahler Ricci soliton is a solution of (\ref{equation-kr_flow})
that moves by a one parameter family of biholomorphisms (which are
induced by a holomorphic vector field). If this vector field comes
from a gradient of a function, we have a gradient K\"ahler Ricci
soliton. In paticular, the equation of a gradient K\"ahler Ricci
soliton $g$ is
\begin{equation}
\label{equation-soliton}
g_{i\bar{j}} - R_{i\bar{j}} = u_{i\bar{j}},
\end{equation}
with $u_{ij} = u_{\bar{i}\bar{j}} = 0$.
\end{definition}

Soliton type solutions are very important ones, since they usually
appear as blow up limits of the Ricci flow and therefore understanding
this kind of solutions helps us understand singularities that the flow
can develop either in a finite or infinite time. Solitons are
generalizations of Einstein metrics (notice that if a function $u$ in
\ref{equation-soliton} is constant, we exactly get a K\"ahler-Einstein
metric). Compactness theorem for K\"ahler Einstein metrics has been
considered and proved by Anderson (\cite{anderson1989}), Bando, Kasue
and Nakajima (\cite{bando1989}), and Tian (\cite{tian1990}). They all
proved that if we start with a sequence of K\"ahler-Eintein metrics
which have uniformly bounded $L^{n/2}$ norms of a curvature, uniformly
bounded Ricci curvatures and diameters, and if $\vol_{g_i}M_i \ge V$
for all $i$, then there exists a subsequence $(M_i,g_i)$ and a compact
K\"ahler Einstein orbifold $(M_{\infty},g_{\infty})$ with a finite set
of singularities, so that $(M_i,g_i)\to (M_{\infty},g_{\infty})$
smoothly away from singular points (in the sense of Gromov-Cheeger
convergence). Moreover, for each singular point there is a
neighbourhood and a gauge so that a lifting of a singular metric in
that gauge has a smooth extension over the origin. The extended smooth
metric is also K\"ahler-Einstein metric. It is natural to expect that
something similar holds in the case of K\"ahler Ricci solitons.
 
Our goal in this paper is to prove the compactness result for the
K\"ahler Ricci solitons, that is, we want to prove the following
theorem.

\begin{theorem}
\label{theorem-compactness}
Let $(M_i,g_i)$ be a sequence of K\"ahler Ricci solitons of real
dimension $n \ge 6$, with $c_1(M_i) > 0$,
\begin{equation}
\label{equation-KR_flow}
\frac{d}{dt}g_i(t) = g_i(t) - \ric(g_i(t)) = \partial\bar{\partial}u_i(t),
\end{equation}
with $\nabla_j\nabla_k u_i = \bar{\nabla}_j\bar{\nabla}_ku_i = 0$,
such that
\begin{enumerate}
\item[(a)]
$\int_{M_i}|\rem|^{n/2}dV_{g_i} \le C_1$,
\item[(b)]
$\ric(g_i) \ge - C_2$,
\item[(c)] 
$A \le \mu(g_i,1/2)$,
\end{enumerate}
for some uniform constants $C_1, C_2, A$, inedependent of
$i$. Then there exists a subsequence $(M_i,g_i)$ converging to
$(Y,\bar{g})$, where $Y$ is an orbifold with finitely many isolated
singularities and $\bar{g}$ is a K\"ahler Ricci soliton in an orbifold
sense (see a definition below).
\end{theorem}

We have a similar result in the case $n=4$, but because of a different 
treatment we will state it separately.

\begin{theorem}
\label{theorem-compactness-4}
Let $(M_i,g_i)$ be a sequence of complex $2$-dimensional K\"ahler
Ricci solitons with $c_1(M_i) > 0$, satisfying
\begin{equation}
\label{equation-KR_flow}
\frac{d}{dt}g_i(t) = g_i(t) - \ric(g_i(t)) = \partial\bar{\partial}u_i(t),
\end{equation}
with $\nabla_j\nabla_k u_i = \bar{\nabla}_j\bar{\nabla}_ku_i = 0$,
such that
\begin{enumerate}
\item[(a)]
$|\ric(g_i)| \le C_1$,
\item[(b)] 
$A \le \mu(g_i,1/2)$,
\end{enumerate}
for some uniform constants $C_1, A$, inedependent of $i$. Then there
exists a subsequence $(M_i,g_i)$ converging to $(Y,\bar{g})$, where
$Y$ is an orbifold with finitely many isolated singularities and
$\bar{g}$ is a K\"ahler Ricci soliton in an orbifold sense.
\end{theorem}

First of all notice that in Theorem \ref{theorem-compactness-4} we
have one condition less. That is because in the case $n=4$ the
condition (a) of Theorem \ref{theorem-compactness} is automatically
fulfilled (a consequence of a Gauss-Bonett formula for surfaces). The
approach that we will use to prove Theorem \ref{theorem-compactness}
is based on Sibner's idea for treating the isolated singularities for
the Yang-Mills equations which requires $n > 4$.  To prove Theorem
\ref{theorem-compactness-4}, we will use the techniques developed by
Uhlenbeck in \cite{uhlen1982} to treat the Yang-Mills equation, and then
later used by Tian in \cite{tian1990} to deal with an Einstein
equation.

\begin{definition}
\label{definition-convergence}
Let $(M_i,g_i)$ be a sequence of K\"ahler manifolds, of real dimension
$n$ ($n$ is taken to be even). We will say that $(M_i,g_i)$ converge
to an orbifold $(M_{\infty},g_{\infty})$ with finitely many isloated
singularities $p_1,\dots,p_N$, where $g_{\infty}$ is a {\it K\"ahler Ricci soliton
in an orbifold sense}, if
\begin{enumerate}
\item[(a)] 
For any compact subset $K\subset
M_{\infty}\backslash\{p_1,\dots,p_N\}$ there are compact sets $K_i\subset M_i$ and 
diffeomorphisms $\phi_i: K_i\to K$ so that $(\phi_i^{-1})^*g_i$ coverge to $g_{\infty}$
uniformly on $K$ and $\phi_{i*}\circ J_i\circ(\phi_i^{-1})_*$ converge to $J_{\infty}$ 
uniformly on $K$, where $J_i,J_{\infty}$ are the almost complex structures of 
$M_i,M_{\infty}$, respectively.
\item[(b)] 
For every $p_i$ there is a neighbourhood $U_i$ of $p_i$ in
$M_{\infty}$ that is covered by a ball $\Delta_r$ in
$\mathrm{C}^{n/2}$ with the covering group isomorphic to a finite
group in $U(2)$. Moreover, if $\pi_i:\Delta_r\to U_i$ is the
covering map, there is a diffeomorphism $\psi$ of $\Delta_r$ so that
$\phi^*\pi_i^*g_{\infty}$ smoothly extends to a K\"ahler Ricci soliton
$C^{\infty}$-metric on $\Delta_r$ in $\mathrm{C}^{n/2}$ with respect
to the standard complex structure.
\end{enumerate}
We will call $(M_{\infty},g_{\infty})$ a {\it generalized K\"ahler
Ricci soliton}.
\end{definition}

We will call $(M_{\infty},g_{\infty})$ a {\it generalized K\"ahler
Ricci soliton}.

The outline of the proof of Theorem \ref{theorem-compactness} is as
follows.

\begin{enumerate}
\item
Obtaining the $\epsilon$-regularity lemma for K\"ahler Ricci solitons
(the analogue of the existing one for Einstein metrics) which says
that a smallness of the $L^{n/2}$ norm of a curvature implies a
pointwise bound on the curvature.
\item
Combining the previous step together with a uniform $L^{n/2}$ bound on
the cuvatures of solitons in our sequence yields a convergence of a
subsequence of our solitons to a toplogical orbifold
$(M_{\infty},g_{\infty})$ with finitely many isolated
singularities. The metric $g_{\infty}$ satisfies the K\"ahler Ricci
soliton equation away from singular points.
\item
Using Moser iteration argument (as in \cite{bando1989} and
\cite{tian1990}) we can show the uniform boundness of
$|\rem(g_{\infty})|$ on $M_{\infty}\backslash\{$ singular points $\}$.
\item
By the similar arguments as in \cite{anderson1989}, \cite{bando1989}
and \cite{tian1990} we can show that $g_{\infty}$ extends to a $C^0$
orbifold metric on $M_{\infty}$.
\item
Using that Ricci potentials of metrics $g_i$ in our sequence are the
minimizers of Perelman's functional $\mathcal{W}$, henceforth
satisfying the elliptic equation, and using harmonic coordinates
around the orbifold points we can show that $g_{\infty}$ extends to a
$C^{\infty}$ orbifold metric on $M_{\infty}$. A lifting of
$g_{\infty}$ above orbifold points is a smooth metric, satisfying a
K\"ahler ricci soliton equation in the covering space.
\end{enumerate}

Due to Perelman, instead of assuming uniform bounds on diameters and
volume noncollapsing constant it is enough to assume condition (c) in
Theorem \ref{theorem-compactness} (see the next section for more
details). In \cite{perelman2002} Perelman has introduced Perelman's
functional 
$$\mathcal{W}(g,f,\tau) = (4\pi\tau)^{-n/2}\int_Me^{-f}[\tau(R|\nabla
f|^2) + f - n]dV_g,$$
under the constraint 
\begin{equation}
\label{equation-constraint}
(4\pi\tau)^{-n/2}\int_Me^{-f}dV_g = 1.
\end{equation} 
He also defined $\mu(g,\tau) = \inf \mathcal{W}(g,\cdot,\tau)$, where
$\inf$ is taken over all functions satisfying the constraint
(\ref{equation-constraint}).

The compactness theorem for K\"aler Einstein manifolds has been
established in \cite{anderson1989}, \cite{bando1989} and
\cite{tian1990}. Almost the same proof of Theorem
\ref{theorem-compactness} yields a generalization of very well known
compactness theorem in K\"ahler Einstein case.

\begin{theorem}
\label{theorem-einstein-cond}
Let $(M_i,g_i)$ be a sequence of K\"ahler Ricci solitons with
$c_1(M_i) > 0$, such that the following holds
\begin{enumerate}
\item[(a)]
$\int_{M_i}|\rem|^{n/2}dV_{g_i} \le C_1$,
\item[(b)]
$\ric(g_i) \ge - C_2$ for $n \ge 6$ and $|\ric(g_i)| \le C_2$ for $n=4$,
\item[(c)] 
$\vol_{g_i}(B_{g_i}(x,r)) \ge \kappa r^n$,
\item[(d)]
$\diam(M_i,g_i) \le C_3$,
\end{enumerate}
for some uniform constants $C_1, C_2, \kappa, C_3$, inedependent of
$i$. Then there exists a subsequence $(M_i,g_i)$ converging to
$(M_{\infty}, g_{\infty})$, where $M_{\infty}$ is an orbifold with
finitely many isolated singularities and $g_{\infty}$ is a K\"ahler Ricci
soliton in an orbifold sense.
\end{theorem}

Due to Perelman's results for the K\"ahler Ricci flow (see
\cite{perelman_kr}), and due to the observation of Klaus Ecker (that
will be discussed in the following section), that was communicated to
the second author by Berhnard List we have the following corollary.

\begin{corollary}
\label{corollary-equivalence}
Theorems \ref{theorem-compactness} and \ref{theorem-compactness-4} are
equivalent to Theorem \ref{theorem-einstein-cond} for dimensions $n\ge
6$ and $n=4$, respectively.
\end{corollary}

The proof of Theorem \ref{theorem-einstein-cond} is essentially the
same as that of Theorem \ref{theorem-compactness}. The only difference
is that we do not have a uniform lower bound on Perelman's functional
$\mu$, so we have imposed uniform bounds on diameters, volume
noncollapsing and a euclidean volume growth, which are also either
implied or given, in the case we start with a sequence of K\"ahler
Einstein manifolds.

{\bf Acknowledgements:} The authors would like to thank R.Hamilton for
amny useful discussions about the problem considered in the paper. The
second author would like to thank J.Viaclovsky for pointing out some
mistakes in the first draft and giving some suggestions how to fix
them.

\end{section}

\begin{section}{Perelman's functional $\mu(g,1/2)$}

We can normalize our solitons, so that $\vol_{g_i}(M) = 1$. In order
to prove convergence we need some sort of $\epsilon$-regularity lemma
(an analogue of $\epsilon$-regularity lemma for Einstein manifolds,
adopted to the case of K\"ahler Ricci solitons). We will use Moser
iteration argument to get a quadratic curvature decay away from
curvature concentration points. Define \\ 
$\mathcal{S}_n(C_2,A,V) =
\{$ K\"aher Ricci solitons $g$ $|\:\: \mu(g,\tau) \ge A$, for all
$\tau\in (0,1)$ and $\ric(g) \le -C_2\}$. \\ 
K\"ahler Ricci shrinking solitons $g_i$ from Theorem
\ref{theorem-compactness} are in $\mathcal{S}_n$.  Due to Perelman we
know that the scalar curvature of a solution $g(t)$ satisfying only
the first defining condition of $\mathcal{S}_n$ is uniformly bounded
along the flow. Perelman also showed this implies $g(t)$ is
$\kappa$-noncollapsed, where $\kappa = \kappa(A)$. In particular,
these bounds are uniform for all elements in $\mathcal{S}_n$ and
\begin{equation}
\label{equation-scalar_bound} 
|R(g_i)| \le \tilde{C},
\end{equation}
for all $i$. 

\begin{lemma}
If $\frac{d}{dt}g(t) = g(t) - \ric(g(t)) = \partial\bar{\partial}u$ is
a shrinking gradient K\"ahler Ricci soliton, then $u(t)$ is a
minimizer of Perelman's functional $\mathcal{W}$ with respect to
metric $g(t)$.
\end{lemma}

\begin{proof}
Let $f(0)$ be a minimizer of $\mathcal{W}$ with respect to metric
$g(0)$. Let $\phi(t)$ be a $1$-parameter family of biholomorphisms
that come from a holomorphic vector field $\nabla u(t)$, such that
$g(t) = \phi(t)^*g(0)$. Function $f(t) = \phi^*f(0)$ is a minimizer of
$\mathcal{W}$ with respect to metric $g(t)$ since $e^{-f(t)}dV_t = dm
= \const$  and since
$$\mu(g(t),\tau) \le \mathcal{W}(g(t),f(t),\tau) = \mathcal{W}(g(0),f(0),\tau)
= \mu(g(0),\tau) \le \mu(g(t),\tau),$$
where the last inequality comes from Perelman's monotonicity for the Ricci
flow. We have that
$$\mathcal{W}(g(0),f(0),\tau) = \mathcal{W}(g(t),f(t),\tau),$$ 
and therefore $\frac{d}{dt}\mathcal{W}(g(t),f(t),\tau) = 0$. On the other
hand, $e^{-f(t)}dV_t = dm = \const$ and by Perelman's monotonicity
formula
\begin{eqnarray*}
0 &=& \frac{d}{dt}\mathcal{W}(g(t),f(t),\tau) \\ 
&=& (4\pi\tau)^{-n/2}\int_M e^{-f(t)}|R_{i\bar{j}} + f_{i\bar{j}} -
g_{i\bar{j}}/(2\tau)|^2dV_t,
\end{eqnarray*}
which implies $R_{i\bar{j}} + f_{i\bar{j}} - g_{i\bar{j}}/(2\tau) = 0$ on $M$,
that is $\Delta f(t) = n/2 - R = \Delta u(t)$ and since $M$ is compact,
$f(t) = u(t)$ (both functions satisfy the same integral normalization 
condition $\int_M e^{-f(t)}dV_t = \int_M e^{-u(t)}dV_t = (4\pi\tau)^{n/2}$). 
\end{proof}

Take $\tau = 1/2$. We have our sequence of K\"ahler Ricci solitons
$(M,g_i)$ which defines a sequence of K\"ahler Ricci flows
$\{g_i(t)\}$, where $g_i\in \mathcal{S}_n(C_2,A,V)$ and $u_i$ is a
Ricci potential for $g_i$. The previous lemma tells us that every
$u_i$ is a minimizer of $\mathcal{W}(g_i,\cdot,1/2)$ and therefore
satisfies,
$$2\Delta u_i - |\nabla u_i|^2 + R(g_i) + u_i - n = \mu(g_i,1/2),$$
which implies $u_i = |\nabla u_i|^2 + \mu(g_i,1/2) + R(g_i) \ge
-\tilde{C}$, by (\ref{equation-scalar_bound}) and condition (c) in
Theorem \ref{theorem-compactness}.  Since we have a uniform lower
bound on $u_i$, as in \cite{perelman_kr} we have that
$$u_i(y,t) \le C\dist_{it}^2(x_i,y) + C,$$
$$|\nabla u_i| \le C\dist_{it}(x_i,y) + C,$$ 
for a uniform constant
$C$, where $u_i(x_i,t) = \min_{y\in M_i}u_i(y,t)$.  In order to prove
that $|u_i(t)|_{C^1} \le C$ for a uniform constant $C$, it is enough
to show that the diameters of $(M_i,g_i(t))$ are uniformly
bounded. Since we have (a), (b), (c) and since $\vol_{g_i}(M) = 1$ for
all $i$, by the same proof as in \cite{perelman_kr} we can show that
the diameters of $(M_i,g_i(t))$ are indeed uniformly bounded.
Therefore, there are uniform constants $C$ and $\kappa$ such that
for all $i$,
\begin{enumerate}
\item
$|u_i|_{C^1} \le C$,
\item
$\diam(M_i, g_i) \le C$,
\item
$|R(g_i)| \le C$,
\item
$(M_i,g_i)$ is $\kappa$-noncollapsed.
\end{enumerate}
This immediatelly implies that Theorems \ref{theorem-compactness} and
\ref{theorem-compactness_4} imply Theorem \ref{theorem-einstein-cond} for
dimensions $n\ge 6$ and $n=4$, respectively. 

A uniform lower bound on
Ricci curvatures, a uniform volume noncollapsing condition and a
uniform upper bound on diameters give us a uniform upper bound on
Sobolev constants of $(M_i,g_i)$, that is, there is a uniform constant
$S$ so that for every $i$ and for every Lipshitz function $f$ on
$M_i$,
\begin{equation}
\label{equation-sobolev}
\{\int_M(f\eta)^{\frac{2n}{n-2}}dV_{g_i}\}^{\frac{n-2}{n}} \le
S\int_M|\nabla(\eta f)|^2dV_{g_i},
\end{equation}
where $\eta$ is a cut off function on $M$.

A uniform lower bound on Ricci curvatures implies the existence of a
uniform constant $V$ such that 
\begin{equation}
\label{equation-volume-up}
\vol_{g_i}(B_{g_i}(p,r)) \le Vr^n,
\end{equation} 
for all $i$, $p\in M$ and all $r > 0$. By Bishop-Gromov volume comparison 
principle we have 
$$\vol_{g_i}(B_{g_i}(p,r)) \le
V_{-C_2}(r)\frac{\vol_{g_i}(B_{g_i}(p,\delta))}{V_{-C_2}(\delta)},$$
which by letting $\delta\to 0$ yields
$$\vol_{g_i}(B_{g_i}(p,r)) \le w_nV_{-C_2}(r) = r^nw_nV_{-C_2r}(1),$$
where $w_n$ is a volume of a euclidean unit ball and $V_{-C_2}(r)$ is
a volume of a ball of radius $r$ in a simply connected space of
constant sectional curvature $-C_2$. The term on the right hand side
of the previous estimate is bounded by $Vr^n$, for a uniform constant
$V$, since $\diam(M,g_i) \le D$ and therefore $0 \le C_2 r \le DC_2$.

Ecker's observation that finishes the proof of Corollary
(\ref{corollary-equivalence}) is as follows.

\begin{lemma}[Ecker]
There is a lower bound on $\mu(g,\tau)$ in terms of a Sobolev constant
$C_S$ for $g$, that is,
$$\mu(g,\tau) \ge -C(n)(1 + \ln C_S(g) + \ln\tau) + \tau\inf_M R(g).$$
\end{lemma}

\begin{proof}
Let $f$ be a minimizer for $\mathcal{W}$, and let $u = \phi^2 =
(4\pi\tau)^{-n/2}e^{-f}$. Then,
\begin{equation}
\label{equation-lower-mu}
\mu(g,\tau) \ge (4\pi\tau)^{-n/2}\int_M (4\tau|\nabla\phi|^2 -
\phi^2\ln\phi^2)dV + \tau\inf_M R(g) - c(n)(1+\ln\tau) 
= I + \tau\inf_M R(g) - c(n)(1+\ln\tau).
\end{equation}
Rescale $g_{\tau} = \frac{g}{4\tau}$, $\phi_{\tau} = (4\tau)^{n/2}\phi$. Then,
$$I = \int_M (|\nabla\phi_{\tau}|^2_{\tau} - \phi_{\tau}^2\ln\phi_{\tau}^2)dV_{\tau},$$
with $\int_M \phi_{\tau}^2dV_{\tau} = 1$. Since a usual Sobolev inequality
(with constant $C_S$) implies a logarithmic Sobolev inequality
$$\int_M (|\nabla w|^2 - w^2\ln w^2)dV_g \ge C(n)(1 + \ln C_S(g)),$$
for every $w \ge 0$ and $\int_M w^2dV_g = 1$. Apply the previous
inequality to $w = \phi_{\tau}$ and to $g_{\tau}$ (note that
$C_S(g_{\tau} \le (1 + 2\sqrt{\tau})C_S(g))$, which implies
$$I \ge -C(n)(1 + \ln((1+ 2\sqrt{\tau})C_S(g))).$$
\end{proof}

Conditions that we are imposing in Theorems \ref{theorem-compactness}
and \ref{theorem-compactness-4} are enough to obtain a uniform upper
bound on Sobolev constants $C_S(g_i)$. This together with (\ref{equation-lower-mu})
gives a uniform lower bound on $\mu(g_i,\tau)$.

\end{section}

\begin{section}{$\epsilon$-regularity lemma for K\"ahler Ricci solitons}
\label{section-regularity}

In this section we will establish $\epsilon$-regularity lemma for
K\"ahler Ricci solitons.

By Bochner-Weitzenbock formulas we have
\begin{equation}
\label{equation-bochner}
\Delta|\rem|^2 = -2\langle\Delta \rem,\rem\rangle + 2|\nabla\rem|^2
- \langle Q(\rem),\rem\rangle,
\end{equation}
where $Q(\rem)$ is quadratic in $\rem$. The Laplacian of a 
curvature tensor in the K\"ahler case reduces to
$$\Delta R_{i\bar{j}k\bar{l}} = \nabla_i\nabla_{\bar{l}}R_{\bar{j}k} + 
\nabla_{\bar{j}}\nabla_kR_{i\bar{l}} + S_{i\bar{j}k\bar{l}},$$
where $S(\rem)$ is quadratic in $\rem$. In the case of a soliton metric
$g\in \mathcal{S}_n(C_2,A,V)$ on $M$, that satisfies 
$g_{i\bar{j}} - R_{i\bar{j}} = \partial_i\bar{\partial}_j u$ , 
by commuting the covariant
derivatives, we get 
\begin{eqnarray}
\label{equation-simpler}
\Delta R_{i\bar{j}k\bar{l}} &=& u_{\bar{j}k\bar{l}i} + u_{i\bar{l}k\bar{j}} +
S_{i\bar{j}k\bar{l}} \nonumber \\
&=& u_{\bar{j}\bar{l}ki} + \nabla_i(R_{\bar{j}k\bar{l}m}u_m) +
u_{ik\bar{l}\bar{j}} + \nabla_{\bar{j}}(R_{i\bar{l}k\bar{m}}u_{\bar{m}}) 
+ S_{i\bar{j}k\bar{l}} \nonumber \\
&=& \nabla_i(R_{\bar{j}k\bar{l}m})u_m + \nabla_{\bar{j}}
(R_{i\bar{l}k\bar{m}})u_{\bar{m}} + S_{i\bar{j}k\bar{l}} \nonumber \\
&=& \nabla\rem * \nabla u,
\end{eqnarray} 
where we have effectively used the fact that $u_{ij} =
u_{\bar{i}\bar{j}} = 0$ and $A*B$ denotes any tensor product of two
tensors $A$ and $B$ when we do not need precise expressions.  By using
that $|u|_{C^1} \le C$ on $M$ and identities (\ref{equation-bochner})
and (\ref{equation-simpler}) we get,
$$\Delta|\rem|^2 \ge -C|\nabla\rem||\rem| + 2|\nabla\rem|^2 -
C|\rem|^3.$$
By interpolation inequality we have
\begin{eqnarray*}
\Delta|\rem|^2 &\ge& (2-\theta)|\nabla\rem|^2 - C(\theta)|\rem|^2 -
C|\rem|^3 \\
&\ge& (2-\theta)|\nabla|\rem||^2 - C(\theta)|\rem|^2 - C|\rem|^3.
\end{eqnarray*}
We will see later how small $\theta$ we will take. Also,
$$\Delta|\rem|^2 = 2\Delta|\rem||\rem| + 2|\nabla|\rem||^2,$$
and therefore,
\begin{equation}
\label{equation-curv}
\Delta|\rem||\rem| \ge -\theta/2|\nabla|\rem||^2 - C(\theta)|\rem|^2
- C|\rem|^3.
\end{equation}
Denote by $u = |\rem|$. Then,
\begin{equation}
\label{equation-u1}
u\Delta u \ge -\theta/2|\nabla u|^2 - C(\theta)u^2 - Cu^3.
\end{equation}

We can now prove the $\epsilon$-regularity theorem for shrinking
gradient K\"ahler Ricci solitons.

\begin{theorem}
\label{theorem-regularity}
Let $g$ be a K\"ahler Ricci soliton as above. Then there exist
constants $C = C(n,\kappa)$ and $\epsilon = \epsilon(n,\kappa)$ so
that for $r < (\frac{\epsilon}{V})^{1/n} = r_0$ if
\begin{equation}
\label{equation-cond1}
\int_{B(p,2r)}|\rem|^{n/2}dV_g < \epsilon,
\end{equation} 
then
\begin{equation}
\label{equation-claim}
\sup_{B(p,r/2)}|\rem(g)|(x) \le \frac{C}{r^2}.
\end{equation}
\end{theorem}

\begin{proof}
It simplifies matters if we assume $r=1$. We may assume that, since
(\ref{equation-curv}), condition (\ref{equation-cond1}) and claim
(\ref{equation-claim}) are all scale invariant. We can start with $r$
so that $\vol B_g(p,2r)$ is sufficiently small (it will become clearer
from the proceeding discussion), then rescale our metric so that
$r=1$, $\vol_g B(p,2)$ is small and $\int_{B_g(p,2)}|\rem|^{n/2} <
\epsilon$ (we are using $g$ to denote a rescaled metric as well).

We will first prove $L^{q'}$ bound on $\rem$, for some $q'>n/2$. Let
$S$ be a uniform Sobolev constant, $V$ as in
(\ref{equation-volume-up}).

\begin{lemma}
\label{lemma-L_bound}
There is $q' > n/2$ so that
$\int_{B_g(p,1)}|\rem|^{q'}dV_g \le C$, where $C = C(n,q',V,S)$.
\end{lemma}

\begin{proof}
Let $\phi$ be a nonnegative cut off function that we will choose later
and $q\ge 2$. Multiply (\ref{equation-u1}) by $\phi^2 u^{q-2}$ and
integrate it over $M$.
\begin{eqnarray}
\label{equation-moser_type}
& & \theta/2\int\phi^2|\nabla u|^2 u^{q-2} + C(\theta)\int\phi^2(u^q +
u^{q+1}) \ge \int \phi^2 u^{q-1}(-\Delta u) \nonumber \\
&\ge& 4(q-1)q^{-2}\int|\phi\nabla(u^{q/2})|^2 + 
 4q^{-1}\int\phi u^{q/2}\langle\nabla\phi,\nabla(u^{q/2})\rangle
\end{eqnarray}
We can write $\frac{4}{q^2}\int\phi^2|\nabla u^{q/2}|^2$ for $\int\phi^2|\nabla
u|^2 u^{q-2}$, apply Schwartz and interpolation inequalities to the
second term on the right hand side of (\ref{equation-moser_type}) and get
\begin{equation}
\label{equation-theta}
2(q-1)q^{-2}\int|\phi\nabla u^{q/2}|^2 \le \frac{2\theta}{q^2}\int
\phi^2|\nabla u^{q/2}|^2 + C(\theta)\int\phi^2(u^q + u^{q+1}) 
+ 2(q-1)^{-1}\int|\nabla\phi|^2u^q.
\end{equation}
Choose $0 < \theta < q-1$ (in particular we can choose $\theta =
(q-1)/2$). Then
\begin{equation}
(q-1)q^{-2}\int|\phi\nabla u^{q/2}|^2 \le C\int\phi^2(u^q + u^{q+1})
+ 2(q-1)^{-1}\int|\nabla\phi|^2u^q,
\end{equation}
where $C = C(n)$.  Using
the Sobolev inequality (let $\gamma = \frac{n}{n-2}$) we obtain
\begin{eqnarray*}
\{\int|\phi u^{q/2}|^{2\gamma}\}^{1/{\gamma}} &\le&
\tilde{C}\{q^2(q-1)^{-1}\int\phi^2(u^{q+1} + u^q) + (2q^2(q-1)^{-2}
+ 1)\int|\nabla\phi|^2u^q\} \\
&\le& \tilde{C}q^2(q-1)^{-1}\{(\int_{\supp\phi}u^{n/2})^{2/n} + 
(\int_{\supp\phi}dV)^{2/n}\}\{\int|\phi u^{q/2}|^{2\gamma}\}^{1/{\gamma}} + \\
&+& \tilde{C}(2q^2(q-1)^{-2} + 2)\int|\nabla\phi|^2u^q, 
\end{eqnarray*}
where $\tilde{C} = \tilde{C}(n,S,V)$. Take $\epsilon =
\tilde{C}^{-1}q^{-2}(q-1)/4$ and let $\phi$ be a cut off function with compact
support in $B_g(p,2)$, equal to $1$ on $B_g(p,1)$ and such that
$|\nabla\phi| \le C$. Take $q = n/2$. Then we get the following
estimate, 
$$\{\int|\phi u^{q'}\}^{1/{\gamma}} \le \bar{C}_1,$$
that is,
\begin{equation}
\label{equation-better_bound0}
\int_{B_g(p,1)}u^{q'}dV_g \le \bar{C},
\end{equation}
where $\bar{C} = \bar{C}(C_2,A,V,n)$ and $q' = q\gamma = n\gamma/2 > n/2$. 
\end{proof}

Let $q'$ be as in Lemma \ref{lemma-L_bound}. Take some $\beta > 1$,
let $\phi$ be a fuction with compact support in $B_g(p,1)$, equal to
$1$ on $B_g(p,1/2)$ and choose $\theta = 1/2$ in
(\ref{equation-theta}). Then the estimate becomes
\begin{eqnarray*}
\int|\phi\nabla u^{\beta/2}|^2 &\le&
\tilde{C}(\frac{\beta^2}{2\beta-3})(\int\phi^2u^{\beta} + \int\phi^2u^{q+1}) +
2\frac{\beta^2}{(\beta-1)(2\beta - 3)}\int|\nabla\phi|^2u^{\beta} \\
&\le& \tilde{C}(\frac{\beta^2}{\beta-1})(\int\phi^2u^{\beta} + \int\phi^2u^{q+1}) +
2\frac{\beta^2}{(\beta-1)^2}\int|\nabla\phi|^2u^{\beta} \\
&\le& \tilde{C}_1(\beta + 1)(\int\phi^2u^{\beta} + \int\phi^2 u^{\beta + 1}) + 
C\int|\nabla\phi|^2u^{\beta}.
\end{eqnarray*}
Furthermore, by (\ref{equation-better_bound0}) we have,
$$\int\phi^2u^{\beta+1} \le (\int_{B_g(p,1)}u^{q'})^{1/q'}(\int|\phi
u^{\beta/2}|^{2\gamma'})^{1/\gamma'} \le \bar{C}_2(\int|\phi
u^{\beta/2}|^{2\gamma'})^{1/\gamma'},$$ 
with $\gamma' =
\frac{q'}{q'-1}$.  By interpolation inequality, with $2^* =
\frac{2n}{n-2} > \frac{2q'}{q'-1} > 2$, since $q' > n/2$, we have
$$||\phi u^{\beta/2}||_{L^{\frac{2q'}{q'-1}}} \le 
\eta||\phi u^{\beta/2}||_{L^{2^*}} + C(n,q')\eta^{-\frac{n}{2q'-n}}
||\phi u^{\beta/2}||_{L^2},$$
for any small $\eta > 0$. By Sobolev inequality (\ref{equation-sob_original}), 
$$||\phi u^{\beta/2}||_{L^{\frac{2q'}{q'-1}}} \le S\eta||\nabla(\phi
u^{\beta/2})||_{L^2} + C(n,q')\eta^{-\frac{n}{2q'-n}}||\phi
u^{\beta/2}||_{L^2}.$$ 
All this yields,
\begin{eqnarray*}
\int_M|\nabla(\phi u^{\beta/2})|^2 &\le&
\tilde{C}_1(1+\beta)\int\phi^2u^{\beta} + (C +
1)\int|\nabla\phi|^2u^{\beta} +
\tilde{C}_1\bar{C}_2(\beta+1)\{\int|\phi u^{\beta/2}|^{2\gamma'}\}^{1/\gamma'} \\
&\le& \tilde{C}_1(1+\beta)\int\phi^2u^{\beta} + (C+1)\int|\nabla\phi|^2u^{\beta}
+ \tilde{C}_1(\beta + 1)\bar{C}_2\eta^2||\phi u^{\beta/2}||^2_{L^{2*}} \\
&+& C(n,q')^2\eta^{-\frac{2n}{2q'-n}}||\phi u^{\beta/2}||^2_{L^2} \\
&\le& \tilde{C}_1(1+\beta)\int\phi^2u^{\beta} + (C+1)\int|\nabla\phi|^2u^{\beta}
+ \tilde{C}_1(\beta + 1)\bar{C}_2 S||\nabla(\phi u^{\beta/2})||^2_{L^2} \\
&+& C(n,q')^2\eta^{-\frac{2n}{2q'-n}}||\phi u^{\beta/2}||^2_{L^2}.
\end{eqnarray*}
Choose $\eta^2 = \frac{1}{3\tilde{C}_1(\beta + 1)\bar{C}_2S}$. Then,
\begin{eqnarray*}
\int|\nabla(\phi u^{\beta/2})|^2 &\le& C_3\int|\nabla\phi|^2u^{\beta} +
C_4(1+\beta + (1+\beta)^{\frac{2n}{2q'-n}})\int\phi^2u^{\beta} \\
&\le& C_5(1+\beta)^{\alpha}\int(|\nabla\phi|^2 + \phi^2)u^{\beta},
\end{eqnarray*}
where $\alpha$ is a positive number depending only on $n$ and $q'$.
Sobolev inequality then implies
$$(\int|\phi u^{\beta/2}|^2)^{1/\gamma} \le C(1+\beta)^{\alpha}
\int(|\nabla\phi|^2 + \phi^2)u^{\beta},$$
where $\gamma = \frac{n}{n-2}$ as before. Let $r_1 < r_2 \le r_0$.
Choose the cut off function as follows. Let 
$\phi \in C_0^1(B_{g_{\infty}}(p,r_2))$ with the property that 
$\phi \equiv 1$ in $B_{g_{\infty}}(p,r_1)$ and 
$|\nabla\phi| \le \frac{C}{r_2-r_1}$. Then we obtain,
$$(\int_{B_g(p,r_1)}u^{\gamma\beta})^{1/\gamma} \le
C\frac{(\beta +
1)^{\alpha}}{(r_2-r_1)^2}\int_{B_{g_{\infty}}(p,r_2)}u^{\beta},$$
that is
$$||u||_{L^{\gamma\beta}(B_g(p,r_1))} \le
(C\frac{(\beta+1)^{\alpha}}{(r_2-r_1)^2})^{1/\gamma}
||u||_{L^{\beta}(B_g(p,r_2))}.$$ 
By Moser iteration technique (exactly
as in the proof of Theorem $4.1$ in \cite{lin}) we get
$$\sup_{B_g(p,1/2)}|\rem|(g)(x) \le
C(\int_{B_g(p,1)}|\rem|^{n/2})^{2/n},$$ 
for a uniform constant $C$. Rescale back to the original metric to get 
(\ref{equation-claim}).
\end{proof}

\end{section}

\begin{section}{Topological orbifold structure of a limit}
\label{section-topology}

By using a quadratic curvature decay proved in Theorem
\ref{theorem-regularity}, in this section we will show that we can
extract a subsequence of $(M_i,g_i)$ so that it converges to an orbifold
in a topological sense. This relies on work by Anderson
(\cite{anderson1989}), Bando, Kasue and Nakajima (\cite{bando1989})
and Tian (\cite{tian1990}). Take $\epsilon_0$ to be a small constant
from Theorem \ref{theorem-regularity}. Define
$$D_i^r =
\{x\in M_i |\:\:\:\int_{B_{g_i}(x,2r)}|\rem|^{n/2}dV_{g_i}
< \epsilon_0\},$$
and similarly
$$L_i^r = \{x\in M_i |\:\:\:
\int_{B_{g_i}(x,2r)}|\rem|^{n/2}dV_{g_i} \ge \epsilon_0\}.$$ 
For each $i$ we can find a maximal $r/2$ separated set,
$\{x_k^i\} \in M_i$, so that the geodesic balls
$B_{g_i}(x_k^i,r/4)$ are disjoint and $B_{g_i}(x_k^i,r)$ form a
cover of $M_i$.  There is a uniform bound on the number of balls $m_i^r$, 
(centred at $x_k^i$, with radius $r$) in
$L_i^r$, independent of $i$ and $r$, which follows from
$$m_i^r\epsilon_0 \le \sum_{k=1}^{m_i(r)}\int_{B_{g_i}(x_k^i,2r)}
|\rem|^{n/2}dV_{g_i} \le m\int_M|\rem|^{n/2}dV_{t_i+t} \le Cm,$$
where $m$ is the maximal number of disjoint balls of radius $r/4$ in
$M_i$ conatined in a ball of radius $4r$, given by
$$m\kappa (r/4)^{n} \le \sum_{k=1}^{m}\vol_{g_i}B_{g_i}(x_k,r/4) \le
\vol_{g_i}B_{g_i}(x,4r) \le Cr^n.$$
By Theorem \ref{theorem-regularity} we have that for all $x\in D_i^r$
and $r \le r_0$,
\begin{equation}
\label{equation-curv_bound00}
|\rem(g_i)|(x) \le \frac{C}{r^2},
\end{equation}
for a uniform constant $C$. This gives the curvatures of $g_i$ being 
uniformly bounded on $D_i^r$, which together with volume noncollapsing condition implies a
uniform lower bound on injectivity radii. We have seen above there is a uniform
upper bound on the number $N$ of points in $(M_i,g_i)$ at which $L^{n/2}$ norm of 
the curvature concentrates. Assume without loss of generality that $N=1$. 
This enables us to assume that $D_i^r = M\backslash B_{g_i}(x_i,2r)$. 
Since K\"ahler Ricci solitons are the solutions of (\ref{equation-KR_flow})
as well, Shi's curvature estimates do apply and therefore by (\ref{equation-curv_bound00}),
$$\sup_{M_i\backslash B_{g_i}(x_i,3r)}|D^k\rem(g_i)| \le C(r,k).$$
Denote by $G_i^r = M_i\backslash B_{g_i}(x_i,3r)$.
By Cheeger-Gromov convergence theorem, we can extract a subsequence so that
$(G_i^r,g_i)$ converges smoothly (uniformly on compact subsets) to  a
smooth open K\"ahler manifold $G^r$  with a metric $g^r$ that satisfies a K\"ahler
Ricci soliton equation $g^r - \ric(g^r) = \partial\bar{\partial}u^r$.   

We now choose a sequence $\{r_j\}\to 0$ with $r_{j+1} < r_j/2$ and
perform the above construction for every $j$. If we set $D_i(r_l) =
\{x\in M|x\in D_i^{r_j}$, for some $j\le l\}$ then we
have
$$D_i(r_l)\subset D_i(r_{l+1}) \subset\dots \subset M_i.$$
For each fixed $r_l$, by the same arguments as above, each sequence
$\{D_i(r_l),g_i\}$ has a smoothly convergent subsequence to a smooth
limit $D(r_l)$ with a metric $g^{r_l}$, satisfying a K\"ahler Ricci
soliton condition.  We can now set $D = \cup_{l=1}^{\infty}D(r_l)$
with the induced metric $g_{\infty}$ that coincides with $g^{r_l}$ on
$D(r_l)$ and which is smooth on $D$.

Following section $5$ in \cite{anderson1989} we can show there are
finitely many points $\{p_i\}$ so that $M_{\infty} = D\cup \{p_i\}$ is
a complete length space with a length function $g_{\infty}$, which
restricts to a K\"ahler Ricci soliton on $D$ satisfying
$$g_{\infty} - \ric(g_{\infty}) = \partial\bar{\partial}u_{\infty},$$
for a Ricci potential $u_{\infty}$ which is a $C^{\infty}$ limit
of Ricci potentials $u_i$ away from singular points. 

To finish the proof of Theorem \ref{theorem-compactness} we still
need to show few things:
\begin{enumerate}
\item [(a)]
There is a finite set of points $\{p_1,\dots,p_N\}$, such that
$$M_{\infty} = D\cup\{p_i\},$$
is a complete orbifold with isolated singularities $\{p_1,\dots,p_N\}$.
\item[(b)] 
A limit metric $g_{\infty}$ on $D$ can be extended to an
orbifold metric on $M_{\infty}$ (denote this extension by $g_{\infty}$
as well). More precisely, in an orbifold lifting around singular
points, in an appropriate gauge, a K\"ahler Ricci soliton equation of
$g_{\infty}$ can be smoothly extended over the origin in a ball in
$\mathrm{C}^{n/2}$.
\end{enumerate}

We will call points $\{p_i\}_{i=1}^N$ {\it curvature singularities} of
$M_{\infty}$ as in \cite{anderson1989}. We want to examine the
structure, topological and metric, of $M_{\infty}$.

The proof that $M_{\infty}$ has a topological structure of an orbifold
is the same as that of \cite{anderson1989}, \cite{bando1989} and
\cite{tian1990} in the case of taking a limit of a sequence of
Einstein metrics, so we will just briefly outline the main
points. Without loosing a generality, assume there is only one
singular point, call it $p$ and assume it comes from curvature
concentration points $x_i\in M$. By a covering argument we may assume
that each ball $B_{g_i}(x_i,2r)$ contains $L_i^r$. By Theorem
\ref{theorem-regularity} we have that
\begin{equation}
\label{equation-before_lim} 
\sup_{M\backslash B_{g_i}(x_i,r)}|\rem(g_i)|(x) \le
\frac{C}{r^2}(\int_{B_{g_i}(x_i,2r)}|\rem|^{n/2}dV_{g_i})^{2/n},
\end{equation}
which by taking limit on $i$ and using a smooth convergence away from 
a singular point $p$ yields,
\begin{equation}
\label{equation-pointwise_bound}
\sup_{M_{\infty}\backslash \{p\}}|\rem|(g_{\infty})(x) \le
\frac{C}{r(x)^2}(\int_{B_{g_{\infty}}(p,2r(x))}|\rem|^{n/2}dV_{g_{\infty}})^{2/n},
\end{equation}
where $r(x) = \dist_{g_{\infty}}(x,p)$. Let $E(r) = \{x\in
M_{\infty}\backslash\{p\} | \:\:\: r(x) \le r\}$. Given a sequence
$s_i\to 0$, let $A(s_i/2,s_i) = \{x\in M_{\infty}|\:\: s_i/2 \le r(x)
\le 2s_i\}$. Rescale the metric $g_{\infty}$ by $s_i^{-2}$. Then the
rescaled Riemannian manifolds $(A(1/2,1),g_{\infty}s_i^{-2})$ have
sectional curvatures converging to zero by
(\ref{equation-pointwise_bound}). There are uniform bounds on the
covariant derivatives $|D^k\rem|$ of the curvature of metrics
$g_{\infty}s_j^{-2}$ on $A(1/2,2)$ for the following reasons:
An estimate (\ref{equation-before_lim}) and Shi's curvature estimates
give us
$$\sup_{M\backslash B_{g_i}(x_i,2r)}|D^k\rem(g_i)|(x) \le
\frac{C}{r^{k+2}},$$
Letting $i\to\infty$ we get
\begin{equation}
\label{equation-derivative_lim}
|D^k\rem|(g_{\infty})(x) \le \frac{C_1}{r(x)^{k+2}},
\end{equation}
for all $x\in M_{\infty}\backslash \{p\}$.
This tells us there are uniform bounds on the covariant derivatives
$|D^k\rem|$ of the curvature of the metric $g_{\infty}s_i^{-2}$ on
$A(1/2,2)$.

As in \cite{anderson1989} and \cite{tian1990} we can get a uiform
bound, independent of $r$, on a number of connected components in
$E(r)$. It follows now that a subsequence
$(A(1/2,2),g_{\infty}s_i^{-2})$ converges smoothly to a flat K\"ahler
manifold $A_{\infty}(1/2,2)$ with a finite number of components. If we
repeat this process for $A(s_i/k,ks_i)$, for any given $k$, passing to
a diagonal subsequence, it gives rise to a flat K\"ahler manifold
$A_{\infty}$. As in \cite{anderson1989}, \cite{bando1989} and
\cite{tian1990} we can show that each component of $A_{\infty}$ is a
cone on a spherical space form $S^{n-1}(1)/\Gamma$ and that every
component is diffeomorphic to $(0,r)\times S^{n-1}/\Gamma$. We will call
$(M_{\infty},g_{\infty})$ a {\it generalized orbifold}.

So far we have proved the following proposition.

\begin{proposition}
Let $(M_i,g_i)$ be a sequence of compact K\"ahler Ricci solitons, with
$c_1(M_i) > 0$, such that $g_i\in \mathcal{S}_n(C_2,A,V)$ and such that
there is a uniform constant $C$,
$$\int_{M_i}|\rem(g_i)|^{n/2}dV_{g_i} \le C.$$ 
There is a subsequence so that $(M_i,g_i)$ converges in the sense of part (a) of definition
\ref{definition-convergence} to a compact generalized orbifold
$(M_{\infty},g_{\infty})$ with finitely many singularities.
Convergence is smooth outside those singular points and $g_{\infty}$
can be extended to a $C^0$ metric in an orbifold sense (in the
corresponding liftings around sigular points).
\end{proposition}

\end{section}

\begin{section}{A smooth metric structure of a limit orbifold $M_{\infty}$
for $n \ge 6$}

In this section we will always assume $n\ge 6$. We will show that a
limit metric $g_{\infty}$ can be extended to an orbifold metric in
$C^{\infty}$ sense. More precisely, let $p$ be a singular point with a
neighbourhood $U\subset M_{\infty}$. Let $U_{\beta}$ be a component of
$U\backslash\{$ singular points$\}$. Recall that each $U_{\beta}$ is
covered by $\Delta_r^* = \Delta_r\backslash \{0\}$. We will show that
in an appropriate gauge, the lifting of $g_{\infty}$ (around singular
points) can be smoothly extended to a smooth metric in a ball
$\Delta_r$ in $\mathrm{C}^{n/2}$. Metric $g_{\infty}$ comes in as a
limit of K\"ahler Ricci solitons $g_i\in \mathcal{S}_n(C_2,A,V)$. A
Sobolev inequality with a uniform Sobolev constant $S$ holds for all
$g_i$. We will show that a Sobolev inequality with the same Sobolev
constant $S$ holds for $g_{\infty}$ as well.

\begin{lemma}
\label{lemma-sobolev_limit}
There is $r_0$ so that for every $r \le r_0$,
\begin{equation}
\label{equation-sob_limit}
(\int_B v^{\frac{2n}{n-2}}dV_{g_{\infty}})^{\frac{n-2}{n}}
\le S\int_B |\nabla v|^2dV_{g_{\infty}},
\end{equation}
for every $v\in C^1_0(B\backslash\{p\})$, where $B = B_{g_{\infty}}(p,r)$.
\end{lemma}

\begin{proof}
Take $r_0$ such that $\vol_{g_{\infty}}B_{g_{\infty}}(p,r_0) \le
Vr_0^n < \epsilon_0$, where $\epsilon_0$ is a small constant from
Theorem \ref{theorem-regularity}. Let $v \in C_{0}^{1}(B\backslash
\{p\})$ and let $\supp(v) = K \subset B\backslash\{p\}$. By the
definition of convergence, there exist diffeomorphisms $\phi_i$ from
the open subsets of $M\backslash \{x_i\}$ to the open subsets of
$M\backslash\{p\}$ that contain $K$, such that every diffeomorphism
$\phi_i$ maps some compact subset $K_{i}$ onto $K$, where $K_{i}$ is
contained in $B_i = B_{g_i}(x_i,r)$, for some sufficiently large $i$
(because of the uniform convergence of metrics on compact subsets). We
have that $\tilde{g_i} =(\phi_{i}^{-1})^*g_i$ converge uniformly and
smoothly on $K$ to $g_{\infty}$.

Let $F_i = \phi_{i}^* (v)$. Then, $\supp F_i \subset K_i \subset
B(x_i,r)\backslash \{x_i\}$. Let $\{\eta_{i}^{k}\}$ be a sequence of
cut-off functions, such that $\eta_{i}^{k} \in C_{0}^{1}(B_i
\backslash \{x_i\})$ and $\eta_{i}^{k} \to 1 (k \to \infty) \ \
\forall i$, and:
$$\int_{B_i}|D\eta_{i}^{k}|^2 \to 0 \ \ (k\to\infty).$$
$\eta_{i}^{k} F_i$ is a function of compact support in
$B_i$. Then by Sobolev inequality:
$$(\int_{B_i}|\eta_{i}^{k} F_i|^{\frac{2n}{n-2}} dV_{g_i})^{\frac{n-2}{n}} \leq
S\int_{B_i}|D(\eta_{i}^{k}F_i)|^2 dV_{g_i}.$$
We can bound $F_i$ with some constant $C_i$ (as a continuous function
on a compact set), and therefore:
$$\int_{B_i}|D(\eta_{i}^{k}F_i)|^2 \leq
S(\int_{B_i}|D\eta_{i}^{k}|^2 C_i +
\int_{B_i}|DF_{i}|^2(\eta_{i}^{k})^2).$$
Let $k$ tend to $\infty$. Then we get:
$$(\int_{{B_i}\backslash \{x_i\}}|F_{i}|^{\frac{2n}{n-2}})^{\frac{n-2}{n}} \leq
C\int_{B_i\backslash \{x_i\}}|DF_i|^2.$$
Since $\supp F_i \subset K_i$, after changing the coordiantes via map
$\phi_i$ we get:
$$(\int_{K}|v|^{\frac{2n}{n-2}} dV_{\tilde{g_i}})^{\frac{n-2}{n}} \leq S \int_{K}|Dv|^2
dV_{\tilde{g_i}}.$$
Metrics $\{\tilde{g_i}\}$ converge uniformly on K to $g_{\infty}$, so if we let
$i$ tend to $\infty$ in the above inequality, keeping in mind that
$\supp v = K\subset B\backslash\{p\}$, we get that:
$$(\int_B|v|^{\frac{2n}{n-2}} dV_{g_{\infty}})^{\frac{n-2}{n}} \leq
S\int_B|Dv|^2 dV_{g_{\infty}}.$$
\end{proof}	

\begin{remark}
\label{remark-sob_ext}
Observe that (\ref{equation-sob_limit}) also holds for $v\in
W^{1,2}(B)$ (by similar arguments as in \cite{bando1989}). Namely,
let $L_k(t)$ be
$$L_k(t) = \left \{
\begin{array}{lll}
k, & \mbox{for } t \ge k, \\
t, & \mbox{for } |t| < k, \\
-k, & \mbox{for } t \le -k.
\end{array}
\right .$$ 
Then,
$$\{\int_B|L_k(v)|^{2\gamma}\}^{1/\gamma} \le S\int_B|\nabla L_k(v)|^2
= S\int_{|v|<k}|\nabla v|^2.$$ 
Letting $k\to\infty$, by Fatou's lemma
we get (\ref{equation-sob_limit}) for $v\in W^{1,2}(B)$
\end{remark}

\begin{subsection}{Curvature bounds in punctured neighbourhoods of singular points}

Choose $r_0' = 2r_0$ as in Lemma \ref{lemma-sobolev_limit}. Decrease
$r_0'$ if necessary so that
$\int_{B_{g_{\infty}}(p,2r_0)}|\rem|^{n/2}dV_{g_{\infty}} < \epsilon$,
where $\epsilon$ is chosen to be small. Since $g_{\infty}$ is a limit
metric of a sequence of K\"ahler Ricci solitons whose curvatures
satisfy (\ref{equation-curv}), we get that $\rem = \rem(g_{\infty})$
also satisfies
\begin{equation}
\label{equation-lim_curv}
\Delta|\rem||\rem| \ge -\theta/2|\nabla|\rem||^2 - C(\theta)|\rem|^2
- C|\rem|^3,
\end{equation} 
for small $\theta\in (0,1)$. Our goal is to show that
the curvature of $g_{\infty}$ is uniformly bounded on
$B_{g_{\infty}}(p,2r_0)\backslash\{p\}$. We will use this
curvature bound to show a smooth extension of a lifting of an
orbifold metric over the origin in $\mathrm{C}^n$. Denote by $u = |\rem|$.
Function $u$ is then a nonnegative solution satisfying
\begin{equation}
\label{equation-u}
u\Delta u \ge -\theta/2|\nabla u|^2 - Cu^2 - Cu^3.
\end{equation}
This is a special case of more general inequality
\begin{equation}
\label{equation-uf}
u\Delta u \ge -\theta/2|\nabla u|^2 - Cfu^2,
\end{equation}
where $f\in L^{n/2}$. By Fatou's lemma we also have that
$$\int_{M_{\infty}}|\rem|^{n/2}dV_{g_{\infty}} \le
\liminf_{i\to\infty}\int_M|\rem|^{n/2}dV_{g_i} \le C.$$

Remember that we are treating the case $n\ge 6$ which allows us to
adopt the approach of Sibner in \cite{sibner1985}. Our goal is to
show that $u\in L^p(B)$, for some $p > n/2$, because it will give us a
uniform bound on the curvature of $g_{\infty}$ away from a singular
point $p$. The proof of the following lemma is similar to the proof of
Lemma $2.1$ in \cite{sibner1985}.

\begin{lemma}
\label{lemma-L_p_punctured}
Let $u \ge 0$ be $C^{\infty}$ in $M_{\infty}\backslash \{p\}$ and
satisfy there (\ref{equation-u}), with $u\in L^{n/2}$. If $u\in
L^{\frac{2nq_0}{n-2}}\cap L^{2q}$, then $\nabla u^q\in L^2$ and in a
sufficiently small ball $B$, for all $\eta\in C^{\infty}_0(B)$,
$$\int_B \eta^2|\nabla u^q|^2 \le C\int_B |\nabla\eta|^2u^{2q}.$$
\end{lemma}

\begin{proof}
As in \cite{sibner1985} we will choose a particulary useful test
function. Let
$$F(u) = \left \{
\begin{array}{ll}
u^q, & \mbox{for }  0 \le u \le l, \\
\frac{1}{q_0}(ql^{q-q_0}u^{q_0} + (q_0-q)l^q) , & \mbox{for } l \le u,
\end{array}
\right .$$ 
and 
$$F_1(u) = \left \{
\begin{array}{ll}
u^{q-1}, & \mbox{for }  0 \le u \le l, \\
\frac{1}{q_0}(ql^{q-q_0}u^{q_0-1} + \frac{(q_0-q)l^q}{u}) , & \mbox{for } l \le u,
\end{array}
\right .$$ 
Set $G(u) = F_1(u)F'(u)$. Then (see \cite{sibner1985}) the
following inequalities are satisfied
\begin{equation}
\label{equation-eq1}
F \le \frac{q}{q_0}l^{q-q_0}u^{q_0},
\end{equation}
\begin{equation}
\label{equation-eq2}
uFF' \le qF^2,
\end{equation}
\begin{equation}
\label{equation-eq3}
(FF')' \ge C'F'^2, \:\:\: C' > 0,
\end{equation}
where the last inequality fails if $q_0 \le 1/2$ (that is the reason
we have assumed $n > 4$ at the moment). 
Let $\eta\in C^{\infty}_0(B)$, for a sufficiently
small ball $B$ and $\bar{\eta} = 0$ in a neighbourhood of $p$. If $\xi$ is
a test function, from (\ref{equation-u}) we have
$$\int\nabla u \nabla (u\xi) \le \theta/2\int|\nabla u|^2\xi + C\int(u^2
+ u^3)\xi.$$ 
In particular, choose $\xi$ to be
$(\eta\bar{\eta})^2F_1(u)F'(u)$. Integrating by parts, 
using (\ref{equation-eq3}), we get
\begin{eqnarray*}
\int|\nabla F(u)|^2(\eta\bar{\eta})^2 - \frac{\theta}{C}\int|\nabla
u|^2(\eta\bar{\eta})|^2 F_1F' &\le& C_1\int\nabla u
FF'(\eta\bar{\eta})\nabla(\eta\bar{\eta}) + \\
&+& C_1\int(u+1)uFF'(\eta\bar{\eta})^2.
\end{eqnarray*}
\begin{eqnarray*}
C_1\int\nabla uFF'(\eta\bar{\eta})\nabla(\eta\bar{\eta}) &=&
C_1\int \nabla F(u) F (\eta\bar{\eta})\nabla(\eta\bar{\eta}) \\
&\le& 1/4\int|\nabla F(u)|^2(\eta\bar{\eta})^2 + 
C_2\int F^2|\nabla(\eta\bar{\eta})|^2.
\end{eqnarray*}
Using Lemma \ref{lemma-sobolev_limit} and (\ref{equation-eq2}) we get,
\begin{eqnarray*}
C_1\int(u+1)uFF'(\eta\bar{\eta})^2 &\le& C_1q\int (u+1)F^2(\eta\bar{\eta})^2 \\
&\le& C_1q\{\int(u+1)^{n/2}\}^{2/n}\{\int (F\eta\bar{\eta})^{\frac{2n}{n-2}}\}^
{\frac{n}{n-2}} \\
&\le& C_1qS||u+1||_{L^{n/2}(B)}||\nabla(\eta\bar{\eta})||_2^2,
\end{eqnarray*}
where $S$ is a Sobolev constant. Since $u+1\in L^{n/2}$, we can choose $B$ small
so that $C_1qS||u+1||_{L^{n/2}(B)} < 1/4$. Then,
\begin{equation}
\label{equation-k_inf}
\int|\nabla F(u)|^2(\eta\bar{\eta})^2 -
\frac{\theta}{C_3}\int|\nabla u|^2(\eta\bar{\eta})^2F_1F' \le C_4\int
F^2|\nabla(\eta\bar{\eta})|^2.
\end{equation}
Choose a sequence $\eta_k\to 1$ on $B$ with $\int|\nabla\eta_k|^n\to 0$ 
as $k\to\infty$. The term we have to estimate is 
\begin{eqnarray*}
\int\eta^2|\nabla\eta_k|^2F^2 &\le& C(l)\int|\nabla\eta_k|^2u^{2q_0} \\
&\le& C(l) \{\int|\nabla\eta_k|^n\}^{2/n}\{\int u^{\frac{2nq_0}{n-2}}\}^{(n-2)/n},
\end{eqnarray*}
which tends to zero as $k\to\infty$, since the last factor on the right is bounded.
If we let $k\to\infty$ in (\ref{equation-k_inf}), we get
\begin{equation}
\int|\nabla F(u)|^2\eta^2 - \frac{\theta}{C_3}\int|\nabla u|^2\eta^2F_1 F'
\le C_4\int|\nabla\eta|^2F^2.
\end{equation}
We will see later we may assume $q_0 = 1$ and $q_0 \le q$. 
Choose small $\theta$ so that:
\begin{itemize}
\item
for $u\le l$,
$$|\nabla F(u)|^2 - \frac{\theta}{C_3}F_1F' |\nabla u|^2 = F'|\nabla
u|^2(F' - \frac{\theta}{C_3}F_1) = u^{q-1}F'|\nabla
u|^2(q-\frac{\theta}{C_3}) \ge 0,$$ 
\item
and for $u\ge l$,
$$|\nabla F(u)|^2 - \frac{\theta}{C_3}F_1F' |\nabla u|^2 = F'|\nabla
u|^2\{qu^{q_0-1}l^{q-q_0}(1-\frac{\theta}{C_3q_0}) -
\frac{(q_0-q)l^q}{u}\} \ge 0.$$ 
\end{itemize}
This implies
\begin{equation}
\label{equation-limit_l}
\int_{u\le l}|\nabla F(u)|^2\eta^2(1 - \frac{\theta}{C_3q}) \le
C_4\int|\nabla\eta|^2F^2.
\end{equation}
For every $l$ we define $F$. Since for $u\ge l$, we have that
$\frac{1}{q_0}(ql^{q-q_0}u^{q_0} + (q_0-q)l^q) \le \frac{q}{q_0}u^q$ and since
$u\in L^{2q}$, for every $\epsilon > 0$ there is $\delta$ so that
whenever $\vol(E) < \delta$, for every $l$ we have 
$\int_E |\nabla\eta|^2F^2 < \epsilon$. Moreover, there is $l_0$
so that for all $l\ge l_0$, we have 
$\vol(\{u\ge l\}) \le \frac{\int_B u^{2q}}{l^q} < \delta$,
which implies,
$$\int|\nabla\eta|^2F^2 = \int_{u\le l}|\nabla\eta|^2F^2 +
\int_{u\ge l}|\nabla\eta|^2F^2 < \int_{u\le l}|\nabla\eta|^2 u^{2q} + \epsilon.$$
Since $F(u) \to u^q$ as $l\to\infty$, letting $l\to\infty$ in 
(\ref{equation-limit_l}) we get
$$\int\eta^2|\nabla u^q|^2 \le C_5\int|\nabla\eta|^2u^{2q} + \epsilon.$$
Since $\epsilon > 0$ can be arbitrarily small, we get
$$\int\eta^2|\nabla u^q|^2 \le C_5\int|\nabla\eta|^2u^{2q}.$$
\end{proof}

\begin{lemma}
Let $u$ be a nonnegative function as above. Then $u\in L^p$,
for some $p > \frac{n}{2}$.
\end{lemma}

\begin{proof}
Since $u = |\rem(g_{\infty})|\in L^{n/2}$ and $n\ge 6$ (we have
assumed the real dimension $n > 4$), we can choose $q_0 = 1$ and
$q=\frac{n}{4}$. Since $u$ is a nonnegative solution of
(\ref{equation-u}), applying Lemma \ref{lemma-L_p_punctured} to $u$,
we find that $\nabla u^{n/4} \in L^2(B)$. By Remark
\ref{remark-sob_ext}, we can apply Sobolev inequality to $u^{n/4}$ to
conclude that $u\in L^p$ with $p = \frac{n}{2}(\frac{n}{n-2}) >
\frac{n}{2}$.
\end{proof} 

Since $\vol_{g_{\infty}}(M_{\infty}) < \infty$, by the previous lemma,
$u\in L^p(B)$, for $p\in (0, \frac{n}{4}\frac{2n}{n-2}]$. Take
$q_0=1$, $q\in [n/2,\frac{n}{4}\frac{2n}{n-2}]$ and repeat the proof
of Lemma \ref{lemma-L_p_punctured} to get $\nabla u^q \in L^2(B)$ for
all such $q$. By Remark \ref{remark-sob_ext} we have $u\in L^s(B)$,
for $s\in [\frac{n}{4}\frac{2n}{n-2}, \frac{n}{4}(\frac{2n}{n-2})^2]$.
If we keep on repeating this, at the $k$-th step we get $\nabla u^q \in L^2$
for $q\in (0, \frac{n}{4}(\frac{2n}{n-2})^k]$ and
$u\in L^q(B)$ for $q\in (0, \frac{n}{4}(\frac{2n}{n-2})^{k+1}]$. Since
$(\frac{2n}{n-2})^k\to\infty$ as $k\to\infty$, we can draw the following 
conclusion.

\begin{lemma}
\label{lemma-all_p}
If we adopt the notation from above, we have $u\in L^q(B)$ and $\nabla
u^q \in L^2(B)$ for all $q$.
\end{lemma}

\begin{remark}
\label{remark-general_f}
We could get the same conclusion for nonnegative functions $u$
satsfying (\ref{equation-uf}) with $f\in L^{n/2}$ (in the case of
$u=|\rem(g_{\infty})|$, $f = u + 1$).
\end{remark}

The previous lemma helps us get the uniform bound on the curvature of
$g_{\infty}$ in a punctured neighbourhood of a singular point $p$, that is,
we have the following proposition in the case $n > 4$.

\begin{proposition}
\label{proposition-curvature_bounds}
There is a uniform constant $C$ and $r_0 > 0$ so that 
$$\sup_{B\backslash\{p\}}|\rem(g_{\infty}|(x) \le \frac{C}{r_0^2},$$
where $B = B_{g_{\infty}}(p,r_0)$.
\end{proposition}

\begin{proof}
Combining Lemma \ref{lemma-all_p} and Remark \ref{remark-sob_ext}, for
any cut off function $\eta$ with compact support in $B$ and any $q$ we
have a Sobolev inequlity 
$$(\int|\eta u^{q/2}|^{2\gamma})^{1/\gamma} \le C \int|\nabla(\eta
u^{q/2})|^2,$$ 
with a uniform constant $C$. The rest of the proof is
the same as the proof of Theorem \ref{theorem-regularity} in the case
of a smooth shrinking K\"ahler Ricci soliton. Choosing $r_0$
sufficiently small, so that $\int_B|\rem|^{n/2}dV_{g_{\infty}}$ and
$\vol_{g_{\infty}}(B)$ are small enough, we get
$$\sup_{B_{g_{\infty}}(p,r/2)\backslash\{p\}}|\rem(g_{\infty}|(x) \le
\frac{C}{r_0^2}.$$
\end{proof}

\end{subsection}

\begin{subsection}{Smoothing property of K\"ahler Ricci solitons}
\label{subsection-same}

In section \ref{section-topology} we have showed that $g_{\infty}$
extends to a $C^0$-metric on $M_{\infty}$ in the sense that each
singular point $p_i\in M_{\infty}$ has a neighbourhood that is covered
by a smooth manifold, diffeomorphic to a punctured ball
$\Delta_r^*\subset \mathrm{C}^{n/2}$. If we denote by $\phi_i$ those
diffeomorphisms and by $\pi_i$ the covering maps, then the pull-back
metric $\phi_i^*\circ\pi_i^*(g_{\infty})$ extends to a $C^0$-metric on
the ball $\Delta_r$. We know that $g_{\infty}$ satisfies a soliton
equation away from orbifold points. Note that the metric
$\phi_i^*\pi^*g_{\infty}$ is a K\"ahler Ricci soliton in $\Delta_r$,
outside the origin. Our goal is to show that $g_{\infty}$ extends to a
$C^{\infty}$-metric on $\Delta_r$. That implies
$\phi_i^*\circ\pi_i^*(g_{\infty})$ satisfies the soliton equation on
$\Delta_r$, that is $g_{\infty}$ is a soliton metric in an {\it
orbifold sense} (see Definition \ref{definition-convergence}).

Using Proposition \ref{proposition-curvature_bounds} in the case $n
\ge 6$ and Proposition \ref{proposition-proposition_remove} in the
case $n = 4$, and harmonic coordinates constructed in \cite{jost1984},
in the same way as in Lemma $4.4$ in \cite{tian1990} and in the proof
of Theorem $5.1$ in \cite{bando1989} we can show that if $r$ is
sufficiently small, there is a diffeomorphism $\psi$ of $\Delta_r^*$
such that $\psi$ extends to a homeomorphism of $\Delta_r$ and
$$\psi^*(g_{\infty})_{i\bar{j}}(x) - \delta_{i\bar{j}} = O(|x|^2),$$
$$\partial_k\psi^*(g_{\infty})_{ij}(x) = O(|x|).$$
This means $\psi^*g_{\infty}$ is of class $C^{1,1}$ on $\Delta_r$,
that is, there are some coordinates in a covering of
a singular point of $M_{\infty}$ in which $g_{\infty}$ extends 
to a $C^{1,1}$-metric (we may assume $g_{\infty}$ is $C^{1,1}$ 
for further consideration). 

\begin{lemma}
Metric $g_{\infty}$ is actually $C^{\infty}$ on $\Delta_r$.
\end{lemma}

\begin{proof}
In section \ref{section-regularity} we have showed that Ricci
potentials $u_i$ satisfy the following equation,
$$2\Delta u_i - |\nabla u_i|^2 + R(g_i) + u_i - n = \mu(g_i,1/2),$$
with $\Delta u_i = n/2 - R(g_i)$ and therefore,
\begin{equation}
\label{equation-let_infty}
\Delta u_i = |\nabla u_i|^2 - u_i + n/2 + \mu(g_i,1/2).
\end{equation}
We have showed that metrics $\{g_i\}$ uniformly and smoothly converge
to a metric $g_{\infty}$  on compact subsets of $M_{\infty}\backslash\{p\}$
(we are still assuming there is only one singular point, a general case
is treated in the same way). The uniform $C^1$ bounds on $u_i$ (see section
\ref{section-regularity}) and a uniform bound on $R(g_i)$, together with
condition $\mu(g_i,1/2) \ge A$ give 
$$A \le \mu(g_i,1/2) \le \tilde{C},$$
for some uniform constant $\tilde{C}$. We can extract a subsequence 
of a sequence of converging metrics $g_i$ so that 
$\lim_{i\to\infty}\mu(g_i,1/2) = \mu_{\infty}$. If we let $i\to\infty$
in (\ref{equation-let_infty}) we get 
\begin{equation}
\label{equation-regularity0}
\Delta u_{\infty} = |\nabla u_{\infty}|^2 - u_{\infty} -
\mu_{\infty} - n/2,
\end{equation}
with $(u_{\infty})_{ij} = 0$ away from a singular point $p$.
Proposition \ref{proposition-curvature_bounds} gives uniform 
bounds on $|\rem(g_{\infty})|$ on $M_{\infty}\backslash\{p\}$ and
therefore,
$$\sup_{M_{\infty}\backslash\{p\}}
|\nabla\bar{\nabla}u_{\infty}|_{g_{\infty}} \le C,$$
for a uniform constant $C$. Since we also have that 
$(u_{\infty})_{ij} = (u_{\infty})_{\bar{i}\bar{j}} = 0$, 
this together with $\sup_{M_{\infty}\backslash\{p\}}
|u_{\infty}|_{C^1(M_{\infty}\backslash\{p\}} \le C$ (which comes
from $|u_i|_{C^1} \le C$) yields
\begin{equation}
\label{equation-C2_bound}
\sup_{M_{\infty}\backslash\{p\}}|u_{\infty}|_{C^2} \le C.
\end{equation}
Since $g_{\infty}$ is $C^{1,1}$ in $\Delta_r$, and since for any two
points $x, y \in \Delta_r^*$ such that a line $\overline{xy}$ does not
contain the origin, due to (\ref{equation-C2_bound}), we have $|\nabla
u_{\infty}(x) - \nabla u_{\infty}(y)| \le C|x - y|$, we can conclude
that $\nabla u_{\infty}$ extends to the origin in $\Delta_r$. Moreover,
$u_{\infty} \in C^{1,1}(\Delta_r)$.

Take the harmonic coordinates $\Phi$ for $g_{\infty}$ (see
\cite{turck1981}) in $\Delta_r$. Outside the origin, $\Phi$ is smooth
and $h = \Phi^*g_{\infty}$ satisfies
\begin{equation}
\label{equation-metric_reg}
\Delta h = 2h - \partial\bar{\partial}u_{\infty}.
\end{equation}

Since $g_{\infty} \in C^1(\Delta_r)$ and $u_{\infty}\in
C^{1,1}(\Delta_r)$, the right hand side of
(\ref{equation-regularity0}) is in $C^{0,1}(\Delta_r)$. By elliptic
regularity this implies $u_{\infty}\in C^{2,\alpha}(\Delta_r)$, for
some $\alpha\in (0,1)$.  We have that $u_{\infty}\in
C^{2,\alpha}(\Delta_r)$ and $g_{\infty}$ is of class
$C^{1,1}(\Delta_r)$ and therefore by results of DeTurck and Kazdan in
\cite{turck1981}, $u_{\infty}$ is at least $C^{1,1}$ and $g_{\infty}$
is of class $C^{1,1}$ in harmonic coordinates in $\Delta_r$. We will
write $g_{\infty}$ for $\Phi^*g_{\infty}$ and from now on when we
mention regularity, or being of class $C^{k,\alpha}$, we will assume
harmonic coordinates. The right hand side of
(\ref{equation-regularity0}) is of class $C^{0,1}(\Delta_r)$, so by
elliptic regularity, $u_{\infty}$ is of class
$C^{2,\alpha}(\Delta_r)$. By elliptic regularity applied to
(\ref{equation-metric_reg}), we get $g_{\infty}$ is
$C^{2,\alpha}(\Delta_r)$. From (\ref{equation-regularity0}) we get
$u_{\infty}$ being of class $C^{3,\alpha}(\Delta_r)$, since the right
hand side of (\ref{equation-regularity0}) is in
$C^{1,\alpha}(\Delta_r)$. Now again from (\ref{equation-metric_reg}),
$g_{\infty}$ is of class $C^{3,\alpha}$ in $\Delta_r$.

If we keep repeating the argument from above, we will obtain that
$g_{\infty}$ is of class $C^{k}$ in $\Delta_r$, for any $k$, that is,
there are coordinates $\psi$ in $\Delta_r$ (disc $\Delta_r$ covers a
neighbourhood of an orbifold point in $M_{\infty}$), such that
$\psi^*\pi^*g_{\infty}$ in those coordinates extends to a
$C^{\infty}$-metric on $\Delta_r$, where $\pi$ is just a covering
map. In particular, the K\"ahler Ricci soliton equation of
$\psi^*\pi^*g_{\infty}$ holds in all $\Delta_r$.
\end{proof}

\end{subsection}

\end{section}

\begin{section}{A smooth orbifold singularity in the case $n=4$}

In this section we deal with the case $n=4$, that is, we want to prove
Theorem \ref{theorem-compactness-4}. The first four sections and the
subsection \label{subsection-same} apply to a four dimensional case, a
different approach has to be taken when one wants to prove the
curvature of a limit metric $g_{\infty}$ is uniformly bounded away
from isolated singularities. We will use the same notation from the
previous sections. To prove that $g_{\infty}$ extends smoothly to a
smooth orbifold metric, satisfying the K\"ahler-Ricci soliton equation
in a lifting around each singular point we will use Uhlenbeck's theory
of removing singularities of Yang-Mills connections in a similar way
Tian used it in \cite{tian1990}.

Assume $M_{\infty}$ has only one singular point $p$. Let $U$ be a
small neighbourhood of $p$ and let $U_{\beta}$ be a component of
$U\backslash\{p\}$. Recall that $U_{\beta}$ is covered by $\Delta_r^*
\subset C^2$ with a covering group $\Gamma_{\beta}$ isomorphic to a
finite group in $U(2)$ and $\pi_{\beta}^*g_{\infty}$ extends to a
$C^0$ metric on the ball $\Delta_r$, where $\pi_{\beta}$ is a covering
map. In order to prove $g_{\infty}$ extends to a smooth metric in a
covering, we first want to prove the boundness of a curvature tensor
$\rem(g_{\infty})$. The proof is similar to that for Yang-Mills
connections in \cite{uhlenbeck1982} with some
modifications. Considerations based on similar analysis can be found
in \cite{tian1990} and \cite{viac2003}. We will just consider
$\Delta_r^*$  as a real $4$-dimensional manifold.

In section \ref{section-topology} we saw there is a gauge $\phi$ so
that by estimates (\ref{equation-pointwise_bound}) and
(\ref{equation-derivative_lim}),
$$||dg_{ij}||_{g_F}(x) \le \frac{\epsilon(r(x))}{r(x)},$$
$$||d(\frac{\partial g_{ij}}{\partial x_k})||_{g_F} \le \frac{\epsilon(r(x))}{r(x)^2},$$
$$||d(\frac{\partial^2 g_{ij}}{\partial x_k\partial x_l})||_{g_F} \le
\frac{\epsilon(r(x))}{r(x)^3},$$ 
in $B_{g_{\infty}}(p,r)$, where $d$ is the exterior
differential on $\mathrm{C}^2 = \mathrm{R}^4$ and $||\cdot||_{g_F}$ 
is the norm on $T^1\mathrm{R}^4$ with respect to the euclidean metric
$g_F$, and $g$ stands for $\phi^*g_{\infty}$.

Let $\tilde{A}$ be a connection form uniquely associated to a metric
$g_{\infty}$ on $\Delta_r^*$, that is $\tilde{D} = d+\tilde{A}$ is the
covariant derivative with respect to $g_{\infty}$. We can view
$\tilde{A}$ as a function in $C^{1,\alpha}(B_0(p,r), so(4)\times
\mathrm{R}^4)$ for $\alpha\in (0,1)$. The following lemma is
essentially due to Uhlenbeck (\cite{uhlen1982}), but the form in which
we will state it below can be found in \cite{tian1990}. It applies to
our case as well, since we also have an $\epsilon$-regularity theorem
as in \cite{tian1990}.

\begin{lemma}
\label{lemma-lemma_choice_gauge}
Let $r$ be sufficiently small. There is a gauge transformation $u$ in
$C^{\infty}(B(p,2r),so(4))$ such that if $D=e^{-u}\tilde{D}e^u = d+A$,
then $d^*A=0$ on $B(p,2r)$, $d^*_{\psi}A_{\psi} = 0$ on $\partial
B(p,2r)$, where $d^*$, $d^*_{\psi}$ are the adjoint operators of the
exterior differentials on $B(p,2r)$ or $\partial B(p,2r)$ with respect
to $g$, respectively. We also have that
$$\sup_{\Delta(r,2r)}(||A||_g(x)) \le \frac{\epsilon(r)}{r},$$
where $\epsilon(r)\to 0$ as $r\to 0$.
\end{lemma}

We have also the following estimates due to Tian (Lemma $4.2$ in
\cite{tian1990}).

\begin{lemma}
\label{lemma-lemma_tian_est}
Let $A$ be the connection given in Lemma
\ref{lemma-lemma_choice_gauge}. For small $r$ we have
$$\sup_{\Delta(r,2r)}||A||_g(x) \le
Cr\sup_{\Delta(r,2r)}||R_{A}||_g(x),$$
$$\int_{\Delta(r,2r)}||A||_g^2(x)dV_g \le Cr^2\int_{\Delta(r,2r)}||R_A||_g^2 dV_g.$$
\end{lemma}

\begin{proposition}
\label{proposition-proposition_remove}
There exist $0< \delta < 1$ and $r > 0$ such that
$$|\rem(g_{\infty})|(x) \le \frac{C}{r(x)^\delta},$$
where $r(x) = \dist_{g_{\infty}}(x,p)$ for $x\in \Delta_r^*$.
\end{proposition}

\begin{proof}
The proof is the same as that of Proposition $4.7$ in
\cite{uhlen1982}. Let $r_1 = r/2, \dots, r_i = \frac{r_{i-1}}{2},
\dots$. Let $A_i$ be a connection given by Lemma
\ref{lemma-lemma_choice_gauge}. As in Lemma $4.3$, in \cite{tian1990},
if we put $\Omega_i = \Delta(r_i,r_{i-1})$, we have
\begin{eqnarray*}
\int_{\Omega_i}||R_{A_i}||^2_{g_{\infty}}dV_{g_{\infty}} &=& 
-\int_{\Omega_i}\langle[A_i,A_i], R_{A_i}\rangle_{g_{\infty}}dV_{g_{\infty}} -
\int_{\Omega_i}\langle A_i, D_i^*R_{A_i}\rangle_{g_{\infty}}dV_{g_{\infty}} \\
&-& \int_{S_i}\langle A_{i\psi}, (R_{A_i})_{r\psi}\rangle_{g_{\infty}}
+ \int_{S_{i-1}}\langle A_{i\psi}, (R_{A_i})_{r\psi}\rangle_{g_{\infty}},
\end{eqnarray*}
where $S_i = \partial \Delta_{r_i}$, $D_i = d + A_i$. If we sum those
identities over $i$, we get
\begin{eqnarray*}
\int_{\Delta(r,2r)}||R(g_{\infty})||^2dV_{g_{\infty}} &=&
- \sum_i \int_{\Omega_i}\langle R_{A_i},[A_i,A_i]\rangle dV_{g_{\infty}}
- \sum_i \int_{\Omega_i} \langle A_i, D_i^*R_{A_i}\rangle dV_{g_{\infty}} \\
&+& \int_{\partial\Delta_r}\langle A_{1\psi}, (R_{A_1})_{r\psi}\rangle.
\end{eqnarray*}
To get the conclusion of Proposition
\ref{proposition-proposition_remove} we proceed in exactly the same
way as in \cite{tian1990}. We do not have a Yang Mills or an Einstein
condition, so we have to use the Ricci soliton equation to estimate a
term $\int_{\Omega_i}\langle A_i,D_i^*R_{A_i}\rangle dV_{g_{\infty}}$
that appears below. We know $g_{\infty}$ satisfies
$$(g_{\infty})_{i\bar{j}} - R_{i\bar{j}} = u_{i\bar{j}},$$ 
with $u_{ij} = u_{\bar{i}\bar{j}} = 0$. By Bianchi identity we have
$D^*\rem = d^{\nabla}\ric$. Since $\nabla_kR_{i\bar{j}} =
u_{i\bar{j},k} = R_{i\bar{j}k\bar{l}}u_{\bar{l}}$, 
and $|\nabla u| \le C$ on $\Delta_r^*$, by Lemma \ref{lemma-lemma_tian_est} 
we have 
\begin{eqnarray*}
\int_{\Omega_i}\langle A_i,D_i^*R_{A_i}\rangle dV_{g_{\infty}} &\le&
(\int_{\Omega_i}||A_i||^2_{g_{\infty}})^{1/2}(\int_{\Omega_i}|d^{\nabla}\ric|^2)^{1/2}
\\
&\le& Cr_i(\int_{\Omega_i}|\rem|^2dV_{g_{\infty}})^{1/2}
(\int_{\Omega_i}|\rem|^2dV_{g_{\infty}})^{1/2} = 
Cr_i\int_{\Omega_i}|\rem|^2dV_{g_{\infty}}.
\end{eqnarray*}
This yields
$$\sum_i\int_{\Omega_i}\langle A_i,D_i^*R_{A_i}\rangle dV_{g_{\infty}} 
\le Cr\int_{\Delta(r,2r)}|\rem|^2dV_{g_{\infty}}.$$
Similarly as in \cite{tian1990} we get
$$|\rem|_{g_{\infty}}(x) \le \frac{C}{r(x)^{\delta}},$$
for $x\in \Delta_r^*$; for sufficiently small $r$ and some $\delta \in (0,1)$.
\end{proof}

Section \ref{subsection-same} applies to the case $n=4$ as well and
that concludes the proof of Theorem \ref{theorem-compactness-4}. 
\end{section}


\begin{thebibliography}{100}

\bibitem{anderson1989} M.Anderson: Ricci curvature bounds and Einstein
metrics on compact manifolds; Journal of the American Mathematical
Society, Volume 2, Number 3 (1989) 455--490.

\bibitem{bando1989} S.Bando, A.Kasue, H.Nakajima: On a construction of
coordinates at infinity on manifolds with fast curvature decay and
maximal volume growth; Invent.math. 97 (1989) 313--349.

\bibitem{jost1984} Jost,J: Harmonic mappings between Riemannian
manifolds; Proceedings of the center for mathematical analysis;
Australian National University, Vol. 4, 1983.

\bibitem{lawson} Bourguignon,J.P, Lawson, H.B.:Stability and isolation
phenomena for Yang-Mills fields; Comm.Math.Phys. 79 (1981) 189--230.
 
\bibitem{cao1985} H.D.Cao: Deformation of K\"ahler metrics to
K\"aher-Einstein metrics on compact K\"ahler manifolds;
Invent. math. 81 (1985) 359--372.

\bibitem{cao2004} H.D.Cao, R.Hamilton, T.Ilmanen: Gaussian densities
and stability for some Ricci solitons;
http://www.math.ethz.ch/~ilmanen/papers/gaussrep.pdf.

\bibitem{turck1981} D.De Turck, J.Kazdan: Some regularity problems in
Riemannian geometry; Ann. Sci. Ec. Norm. Super. 14 (1981) 249--260.
 
\bibitem{glickenstein2002} D. Glickenstein: Precompactness of
solutions to the Ricci flow in the absence of injectivity radius
estimates; preprint arXiv:math.DG/0211

\bibitem{hamilton1995} R. Hamilton: A compactness property for
solutions of the Ricci flow, Amer. J. Math. 117 (1995) 545--572.

\bibitem{hamiltonMA} R. Hamilton: The formation of singularities in
the Ricci flow, Surveys in Differential Geometry, vol. 2,
International Press, Cambridge, MA (1995) 7--136.

\bibitem{hamilton1999} R. Hamilton: Non-singular solutions of the Ricci
flow on 3 manifolds, Communications in Analysis and Geometry vol. 7
(1999) 695--729.

\bibitem{lin} F.Lin: Elliptic Partial Differential Equations; American
Mathematical Society; Courant Institute of Mathematical Sciences;
ISBN 0-82182-691-3.

\bibitem{nakagawa1993} Y. Nakagawa: An isoperimetric inequality for
orbifolds; Osaka J. Math. 30 (1993), 733--739.

\bibitem{perelman2002} G. Perelman: The entropy formula for the Ricci
flow and its geometric applications; arXiv:math.DG/0211159.

\bibitem{perelman_kr}N.Sesum, G.Tian:Bounding scalar curvature and diameter 
along the K\"ahler Ricci flow (after Perelman); notes in preparation.

\bibitem{natasa_kr} N.Sesum: Convergence of a K\"ahler Ricci flow;
preprint, arXiv:math.DG/0402238 (submitted to the Mathematical
Research Letters).

\bibitem{shi1989} W.-X.Shi: Deforming the metric on complete
Rielamnnian manifolds; J.Diff.Geom. 30 (1989) 223--301.

\bibitem{sibner1985} L. M. Sibner: Isolated point singularity problem;
Mathematische Annalen 271 (1985) 125--131.

\bibitem{tian1990} G.Tian: On Calabi's conjecture for complex surfaces
with positive first Chern class; Inventiones math. 101 (1990),
101--172.

\bibitem{viac2003} G.Tian, J.Viaclovsky: Bach-flat asymptotically
locally euclidean metrics; arXiv:math.DG/0310302 v2.

\bibitem{uhlen1982} K.Uhlenbeck: Removable singularities in Yang-Mills
fields; Communications in Mathematical Physics 83 (1982), 11-29.

\bibitem{wei1996} X. Dai, G.Wei, R.Ye: Smoothing Riemannian metrics
with Ricci curvature bounds; Manuscripta Math. 90 (1996), 49--61.

\end{thebibliography}
\end{document}